\documentclass[12pt]{amsart}
\usepackage{geometry} 
\newcounter{dummy} \numberwithin{dummy}{section}
\newtheorem{theorem}[dummy]{Theorem}
\newtheorem{definition}[dummy]{Definition}
\newtheorem{proposition}[dummy]{Proposition}
\newtheorem{lemma}[dummy]{Lemma}
\newtheorem{note}[dummy]{Note}
\usepackage{graphicx}
\usepackage{amssymb}
\usepackage{epstopdf}
\usepackage{lscape}
\usepackage{indentfirst}
\usepackage{latexsym}
\usepackage{multirow}
\usepackage{tabls}
\usepackage{amsthm}
\usepackage{amsthm,amsmath}
\usepackage{graphicx}
\usepackage{lscape}
\usepackage{indentfirst}
\usepackage{latexsym}
\usepackage{multirow}
\usepackage{tabls}
\usepackage{wrapfig}
\usepackage[all]{xy}
\usepackage[numbers]{natbib}
\usepackage{tikz}
\usepackage{mathrsfs}
\newtheorem{corollary}[dummy]{Corollary}

\usepackage[OT2,T1]{fontenc}
\DeclareSymbolFont{cyrletters}{OT2}{wncyr}{m}{n}
\DeclareMathSymbol{\Sha}{\mathalpha}{cyrletters}{"58}

\title[The $3-$part of the ideal class group]{The 3-part of the ideal class group of a certain family of Real Cyclotomic Fields}
\author{eleni agathocleous}

\begin{document}
\begin{center}
\maketitle{\textbf{Abstract}}
\end{center}

\tiny
In this paper we study the structure of the $3-$part of the ideal class group of a certain family of real cyclotomic fields with $3-$class number exactly $9$ and conductor equal to the product of two distinct odd primes. We employ known results from Class Field Theory as well as theoretical and numerical results on real cyclic sextic fields, and we show that the $3-$part of the ideal class group of such cyclotomic fields must be cyclic. We present four examples of fields that fall into our category, namely the fields of conductor $3 \cdot 331$, $7 \cdot 67$, $3 \cdot 643$ and $7 \cdot 257$, and they are the only ones amongst all real cyclotomic fields with conductor $pq \leq 2021$. The $3-$part of the class number for the two fields of conductor $3 \cdot 643$ and $7 \cdot 257$ was up to now unknown and we compute it in this paper.
\\  
 


\normalsize

\section{Introduction}\label{Intro}
 
The goal of this paper is to investigate the structure of the ideal class group of a certain family of real cyclotomic fields $\mathbb{Q}(\zeta_{pq})^{+}$ of conductor $pq$, $p < q$ distinct odd primes, with $3-$class number exactly 9.

We ask that $p$ (or $q$) $\equiv 1~({\rm mod}~3)$ and we study the sextic field $K_{6} = K_{2} K_{3}  \subseteq \mathbb{Q}(\zeta_{pq})^{+}$, where $K_{2}$ is a real quadratic subfield of conductor $f_{2}$ (therefore $f_{2} = pq, p$ or $q$) and $K_{3}$ is the real cyclic cubic subfield of conductor $f_{3} = p$ (or $q$). We ask further that $f_{2} \neq f_{3}$ and that the class number of $K_{2}$ is divisible by $3$, a condition which is very fast and easy to check with PARI/GP.

With specific hypotheses on the units of the sextic field $K_{6} = K_{2} K_{3}  \subseteq \mathbb{Q}(\zeta_{pq})^{+}$, we show that the $3-$part of the ideal class group of $K_{6}$ must be cyclic of order $9$. We then show that the ideal class group of the associated real  cyclotomic field must also be cyclic. 

In Sections~\ref{Notation} and ~\ref{CapitulationKernel} we state some facts from Class Field Theory and some results on Real Cyclic Sextic Fields, necessary for the proof of Proposition~\ref{Prop1}. In Section~\ref{Examples} we combine numerical data on Real Cyclic Sextic Fields from M$\ddot{a}$ki \cite{Maki1} and from Agathocleous \cite{Aga} and we  present four real cyclotomic fields that fall into our category; namely the fields of conductor $pq$ = $3 \cdot 331$, $7 \cdot 67$, $3 \cdot 643$ and $7 \cdot 257$. These are the only ones amongst all real cyclotomic fields with conductor $pq \leq 2021$.

The $3-$part of the class number for the fields of conductor $3 \cdot 643$ and $7 \cdot 257$ was up to now unknown. Applying the first two steps from \cite{Aga} and the numerical data presented in \cite{Maki1}, we are able to show that the $3-$part of their class number is exactly $9$. The numerical data from applying the method from \cite{Aga} are presented in the Appendix.

\section{Notation and Preliminaries} \label{Notation}

\subsection{Some Basic Class Field Theory Results} \label{ClassFieldTheory}
Let $L/K$ be an extension of number fields with ideal class groups $Cl(L)$ and $Cl(K)$ and corresponding $l-$parts $Cl(L)_{l}$ and $Cl(K)_{l}$. Let $\mathfrak{t}$ denote the transfer of ideal classes $\mathfrak{t} : Cl(K) \mapsto Cl(L)$ induced by mapping an ideal $\mathfrak{c}$ to $\mathfrak{c}O_{L}$, where $O_{L}$ the ring of integers of $L$. Then $\mathfrak{t}$ is a homomorphism with kernel $ker (\mathfrak{t})_{K \mapsto L}$. 

\begin{proposition} \label{Lem1} (\cite{Lemm2}, Proposition 10) Let $L/K$ be a ramified cyclic extension of prime degree $l$. Put $r_{K}$ = rank($Cl(K)_{l})$ and $r_{\mathfrak{t}}$ = rank($ker (\mathfrak{t})_{K \mapsto L})$. Then $|Cl(L)_{l}| \geq l^{r_{K} - r_{\mathfrak{t}}} |Cl(K)_{l}|$.   $\square$  \end{proposition} 

\begin{theorem} \label{Lem2} (\cite{Lemm2}, Theorem 6) If $[Cl(L)_{l} : \mathfrak{t}(Cl(K)_{l})] = l^{a}$ for some $a \leq l - 2 + r_{\mathfrak{t}}$, then $\mathfrak{t}(Cl(K)_{l}) = Cl(L)_{l}^{l}$.   $\square$  \end{theorem}

Let $N_{L/K} : Cl(L) \mapsto Cl(K)$ denote the norm map. By Class Field Theory, $N$ is onto when the extension is ramified. In the case that $L/K$ is a ramified extension of degree $l^{a}$, some $a \geq 0$, and their $l$-class numbers are equal, then $N_{L/K} : Cl(L)_{l} \mapsto Cl(K)_{l}$ is an isomorphism and we have the following proposition. 

\begin{proposition} \label{Wash} (\cite{Wash2}, Lemma 4 (i)) Suppose $L/K$ is an extension of number fields with no nontrivial unramified subextensions $M/K$. Assume $[L : K] = l^{a}$ and $l \nmid h(L)_{l} / h(K)_{l}$. Then the kernel of the map $\mathfrak{t} : Cl(K) \mapsto Cl(L)$ defined above, is exactly the classes of order dividing $l^{a}$. $\square$ \end{proposition}

\begin{note} \label{Note} In the case that $l \nmid [L : K]$, the composition of $N_{L/K}$ and $\mathfrak{t}$ gives $N_{L/K} \circ \mathfrak{t} (\mathfrak{c}) = \mathfrak{c}^{\small{[L : K]}}$ hence the map $Cl(K)_{l} \mapsto Cl(K)_{l} : \mathfrak{c} \mapsto \mathfrak{c}^{\small{[L : K]}}$ is an isomorphism. The maps  $\mathfrak{t} : Cl(K)_{l} \mapsto Cl(L)_{l}$ and $N : Cl(L)_{l} \mapsto Cl(K)_{l}$ are therefore an injection and a surjection respectively and $Cl(K)_{l}$ is isomorphic to a direct summand of $Cl(L)_{l}$. $\square$ \end{note}

\subsection{Real Cyclic Sextic Fields} \label{Sextic}

For the theory and facts on real cyclic sextic fields that we present here, we follow the notation in \cite{Maki1} and \cite{Maki2}.

For $n, m \in \{1, 2, 3, 6\}$, we denote by $K_{n}$ a real cyclic extension of degree $n$ over $\mathbb{Q}$, having conductor $f_{n}$, generating character $\chi_{n}$, class number $h_{n}$, ring of integers $O_{n}$ and unit group $U_{n}$. To follow the same notation, we write $Cl_{n}$ for the ideal class group of $K_{n}$ and $Cl_{n,l}$ for its $l$-part. We denote the $l$-part of the class number by $h_{n,l}$. We denote by $N_{n/m}$ the norm from $K_{n}$ to $K_{m}$. $G = \langle \sigma \rangle$ is the Galois group of $K_{6}$, and the conjugates of a number $\gamma \in K_{6}$ are $\gamma, \gamma' = \gamma^{\sigma}, \gamma'' = \gamma^{\sigma^{2}}$, etc. 

The fundamental units of $K_{2}$ and $K_{3}$ are the norm-positive units denoted by $\mu$ and $\tau$ respectively. The group of relative units is defined as $U_{R} = \{\epsilon \in U_{6} | N_{6/3}(\epsilon) = \pm 1, N_{6/2}(\epsilon) = \pm 1 \}$ and we have that $h_{6} = h_{2} h_{3} h_{R}$ whereas $[U_{6} : U_{2} U_{3} U_{R}] =: Q_{K}$ is known as Leopoldt's Unit Index. The positive integer $h_{R}$ is called the relative class number of $K_{6}$ and it equals $h_{R} = [U_{6} : U_{6}^{*}] [U_{6}^{*} :  \langle -1, \mu, \tau, \tau', \xi_{A}, \xi_{A}' \rangle]$, where $U_{6}^{*} = \langle -1, \mu, \tau, \tau', \xi_{A}, \xi_{A}', \xi_{R}, \xi_{R}' \rangle$. The relative unit $\xi_{R} \in K_{6}$ serves as the generating unit for $U_{R}$ and its existence is proved in M$\ddot{a}$ki (\cite{Maki1}, Ch. 4). Following Hasse, M$\ddot{a}$ki (\cite{Maki1}, p.16) defines the cyclotomic unit $\eta$ of $K_{6}$ as the quotient $\eta = \frac{\xi}{\xi'}$, of a special integer $\xi \in \mathbb{Q}(\zeta_{2f_{6}})$ over its conjugate $\xi'$, where $\xi$ is always a unit in our case since $f_{6} = p q$ is not a prime power. We take $\xi_{A} = \xi$ if $\xi \in K_{6}$ or $\xi_{A} = \eta$ otherwise. Since $\xi_{A} \in K_{6}$ we can write $N_{6/3} (\xi_{A}) = \pm \tau^{u} \tau'^{v}$ and $N_{6/2} (\xi_{A}) = \pm \mu^{w}$.

\begin{theorem}
\label{Th1}
(\cite{Maki1}, Theorem 4, p.17) 

Let $u, v, w$ be as above. Then \[ \langle - 1 \rangle N_{6/3}(U_{6}^{*}) = \left\{\begin{array}{r@{\quad,\quad}l} U_{3} & if \ 2 \nmid u \ or \ 2 \nmid v  \\ \langle -1, \tau^{2}, \tau'^{2} \rangle & if \ 2 | u \ and \ 2 | v \\
\end{array}\right.
\]
and  
\[ N_{6/2}(U_{6}^{*}) = \left\{\begin{array}{r@{\quad,\quad}l} U_{2} & if  \  3 \nmid w  \\ \langle -1, \mu^{3} \rangle & if  \  3 | w \\
\end{array}\right.
\] $\square$
\end{theorem}

\begin{theorem}
\label{Th2}
(\cite{Maki1}, Theorem 6, p.18) 

Let $\langle - 1 \rangle N_{6/3}(U_{6}^{*}) = U_{3}$ and $N_{6/2}(U_{6}^{*}) \neq U_{2}$. If either one of the equations $$x^{3} = \mu \xi_{R} \xi'_{R}, \ x^{3} = \mu^{-1} \xi_{R} \xi'_{R} \ \ \ \ \  (\ast)$$ has a solution $x = \xi_{B} \in K_{6}$, then $[U_{6} : U_{6}^{*}] = 3$ and $N_{6/2}(U_{6}) = U_{2}$. $\square$

\end{theorem}

We notice that when the hypotheses of Theorem~\ref{Th2} hold, since  $\langle - 1 \rangle N_{6/3}(U_{6}^{*}) = U_{3}$, we have that $N_{6/3}(\xi_{B}) = \pm 1$ implying that $\langle - 1 \rangle N_{6/3}(U_{6}) = U_{3}$ as well.

\begin{theorem}
\label{Th2(a)}
(\cite{Maki1}, Theorem 8, p.20) 

Let $\langle - 1 \rangle N_{6/3}(U_{6}^{*}) \neq U_{3}$ and $N_{6/2}(U_{6}^{*}) \neq U_{2}$. 

(i) If either one of the equations $(\ast)$ above has a solution $x = \xi_{B} \in K_{6}$, and one of the equations $$x^{2} = |\tau \xi_{R}|, \ x^{2} = |\tau' \xi_{R}|, \ x^{2} = |\tau \tau' \xi_{R}| \ \ \ \ \ \ (\ast \ast)$$ has a solution $x = \xi_{C} \in K_{6}$, then $[U_{6} : U_{6}^{*}] = 12$, $\langle - 1 \rangle N_{6/3}(U_{6}) = U_{3}$ and $N_{6/2}(U_{6}) = U_{2}$. 

(ii) If one of the equations $(\ast)$ above has a solution $x = \xi_{B} \in K_{6}$, and none of the equations $(\ast \ast)$ above has a solution $x = \xi_{C} \in K_{6}$, then $[U_{6} : U_{6}^{*}] = 3$ and $N_{6/2}(U_{6}) = U_{2}$. $\square$

\end{theorem}

\begin{definition}
\label{def1}
The character $\chi_{6}$ is decomposable if $f_{3} \nmid 3 f_{2}$. If $f_{3} | 3 f_{2}$, then $\chi_{6}$ is decomposable if and only if $f_{2}$ is decomposable.

$f_{2}$ is nondecomposable if and only if 

     \begin{align}
       & f_{2} = q_{0} \ ; \ q_{0}\equiv 1~({\rm mod}~4) \\
       & f_{2} = q_{0} q_{1} \ ; \ q_{0}, q_{1}\equiv 3~({\rm mod}~4) \\
       & f_{2} = 4 q_{1} \ ; \  q_{1}\equiv 3~({\rm mod}~4) \\
       & f_{2} = 8 \\
       & f_{2} = 8 q_{1} \ ; \  q_{1}\equiv 3~({\rm mod}~4).
     \end{align}

\end{definition}

Let $k  = 
\begin{cases}
      0, & \text{if} \  \chi_{6} \  \text{is decomposable} \\
      1, & \text{if} \  \chi_{6} \  \text{is nondecomposable}
    \end{cases}$

\begin{theorem}
\label{Th3}
(\cite{Maki2}, Theorem 2, p.597)

Let $f_{6} = 3^{\lambda} p_{1} \cdot \cdot \cdot p_{\nu} f_{2}$ where $\lambda \in \{0, 1, 2 \}$, $\nu \geq 0$, and the $p_{i}$'s are distinct primes $\equiv 1~({\rm mod}~6)$. 

(i) If there exists an index $i \in \{1, 2, ..., i \}$ such that $(f_{2}/p_{i}) = 1$ or if $\lambda = 2$ and $f_{2}\equiv  1~({\rm mod}~3)$, then $w = 0$.   

(ii) Otherwise, $w = -2^{\nu + k + max \{ 1, \lambda \} - 2} h_{2}$  $\square$ \end{theorem}

\begin{theorem}
\label{Th4}
(\cite{Maki2}, Theorem 3, p.597)
Let $f_{6} = p_{1}^{\lambda} p_{2} \cdot \cdot \cdot p_{\nu} f_{3}$, where the $p_{i}$'s are distinct primes, and $\lambda =1$ if $\nu \geq 0$ and $f_{6}$ is odd, while $p_{1} = 2$ and $\lambda \in \{2, 3 \}$ if $f_{6}$ is even.

(i) If there exists an index $i \in \{1, 2, ..., \nu \}$ such that $\chi_{3}(p_{i}) = 1$, then $(u, v) = (0, 0)$. 

(ii) Otherwise $u^{2} - uv + v^{2} = 3^{\nu + k -1} h_{3}$.

  $\square$ \end{theorem}

\section{Classes of $K_{2}$ in $K_{6}$} \label{CapitulationKernel}

We can assume that $N_{6/2}(\xi_{R}) =1$ since the sign of $\xi_{R}$ can be changed. For $n$ = 2 or 3 we let $U_{Rn} = \{ \epsilon \in U_{6} | N_{6/n}(\epsilon) = 1 \}$. We denote the group of ideal classes of $K_{2}$ that become principal in $K_{6}$ by $\mathscr{C}_{6/2} = \{ Cl_{2,3} (\mathfrak{c}) | \mathfrak{c} O_{6} = \gamma O_{6} \}$, where $\mathfrak{c}$ is a nonzero fractional ideal and $\gamma$ is some element of $K_{6}^{*}$.

The map $g_{2}(Cl_{2,3}(\mathfrak{c})) = \gamma^{1 - \sigma^{2}} U_{6}^{1 - \sigma^{2}}$ is a well-defined homomorphism $g_{2} : \mathscr{C}_{6/2} \mapsto U_{R2}/U_{6}^{1 - \sigma^{2}}$. Assume that $\gamma^{1 - \sigma^{2}} = \epsilon^{1 - \sigma^{2}}$, $\epsilon \in U_{6}$.  Then $\frac{\gamma}{\epsilon} = (\frac{\gamma}{\epsilon})^{\sigma^{2}}$ therefore $\frac{\gamma}{\epsilon} \in K_{2}$ and $\mathfrak{c} = (\frac{\gamma}{\epsilon})O_2$ belongs to the principal ideal class, hence $g_{2}$ is injective. Now assume that an element $\epsilon U_{6}^{1 - \sigma^{2}}$ is in the image $Im(g_{2})$. Then by Hilbert's Theorem 90, there is a $\gamma \in O_{6}^{*}$ such that  $\epsilon = \gamma ^{1 - \sigma^{2}}$. But now $\sigma^{2}$ fixes $\gamma$ hence $\gamma = \mathfrak{c} O_{6} \times \mathfrak{p}_{1}^{\mathfrak{\nu}_{1}} \times ... \times \mathfrak{p}_{h}^{\mathfrak{\nu}_{h}}$ where $\mathfrak{c}$ is an ideal of $O_{2}$ and the $\mathfrak{p}_{i}$'s are distinct prime ideals of $O_{6}$ which ramify in $K_{6}/K_{2}$.  We may assume that for each $i, \mathfrak{p}_{i} \nmid \mathfrak{c}$ so that $\mathfrak{\nu}_{i} = \mathfrak{\nu}_{\mathfrak{p}_{i}}(\gamma)$ is the highest power of $\mathfrak{p}_{i}$ dividing $\gamma$. Then it is easy to see that $\epsilon U_{6}^{1 - \sigma^{2}} \in Im(g_{2})$ if and only if $3 | \mathfrak{\nu}_{i}$ for $i = 1, ..., h$. (\cite{Maki2}, p. 598). 

Part of the proof of both Theorem~\ref{Th4} above as well as Lemma~\ref{Lemma1} below, are based on a fact from \cite{Gras} regarding the index of a subgroup of $U_{3}$ in $U_{3}$ itself. For any unit $\pm 1 \neq \gamma \in U_{3}$, denote by $I(\gamma)$ the index $[U_{3} : \langle -1, \gamma, \gamma' \rangle]$. $U_{3}$ is a module over $\mathbb{Z}[\langle \sigma^{2} \rangle]/(1 + \sigma^{2} + \sigma^{4}) \cong \mathbb{Z}[j]$, where $j$ is a cube root of unity. Following the notation in Gras \cite{Gras}, let $\mathbb{Q}' = \mathbb{Q}(\sqrt{-3})$. Then, as a corollary of Proposition 1 in \cite{Gras}, we have that \begin{equation} I(\gamma^{\lambda + \mu \sigma}) = I(\gamma^{\lambda + \mu j}) = N_{\mathbb{Q}'/\mathbb{Q}}(\lambda + \mu j) \cdot I(\gamma). \end{equation} In particular, $N_{\mathbb{Q}'/\mathbb{Q}}(\lambda + \mu j) = (\lambda + \mu j)(\lambda + \mu j') = \lambda^{2} +\lambda \mu (j + j') + \mu^{2}jj' = \lambda^{2} - \lambda \mu + \mu^{2}$.

The following two facts from Hasse (in \cite{Maki2}, p.600) are also essential for Lemma~\ref{Lemma1}:
$$[U_{R2} : U_{6}^{1 - \sigma^{2}}] = 3^{2 - q}$$ $$3^{q} = [N_{6/2}(U_{6}) : U_{2}^{3}] = 3, \  \text{if} \ N_{6/2}(U_{6}) = U_{2}.$$

\begin{lemma} \label{Lemma1} (\cite{Maki2}, Lemma 2, p.601) If $N_{6/2}(U_{6}) = U_{2}$, then $\{ \tau^{i} | i = 0, 1, 2 \}$ is a system of coset representatives of $U_{R2}$ with respect to the subgroup $U_{6}^{1 - \sigma^{2}}$. $\square$ \end{lemma}

Lemma~\ref{Lemma1} is essential for the proof of Theorem~\ref{3Capitulation} below, regarding capitulation in $K_{6}/K_{2}$.

\begin{theorem} \label{3Capitulation} (\cite{Maki2}, Theorem 5, p.601) If $N_{6/2}(U_{6}) = U_{2}$ then $\mathscr{C}_{6/2} = 1$ and there is no capitulation in $K_{6}/K_{2}$. \end{theorem}

\emph{Proof:} From Hilbert's Theorem 90 we have that $\tau^{i} = \gamma^{1 - \sigma^{2}}$ for some $\gamma \in O_{6}^{*}$. If $\tau^{i}U_{6}^{1 - \sigma^{2}} \in Im(g_{2})$, then there is an ideal $\mathfrak{c}$ of $O_{2}$ such that $\mathfrak{c}O_{6} = \gamma O_{6} \Rightarrow \mathfrak{c}^{1 + \sigma^{3}} O_{6} = \gamma^{1 + \sigma^{3}} O_{6}$ since $\sigma^{3}$ is an automorphism. Now for $\mathfrak{c} \in O_{2}$ we have that $\mathfrak{c}^{1 + \sigma^{3}} O_{6} = c O_{6}$ for some $c \in \mathbb{Z}$ since $1 + \sigma^{3} = N_{K_{2}/\mathbb{Q}}$. Hence, there must be an $\epsilon \in U_{3}$ such that $\gamma^{1 + \sigma^{3}} =  c \epsilon$. But this would give the equation $\tau^{2i} = (\tau^{i})^{1 + \sigma^{3}} = \epsilon^{1 - \sigma^{2}}$. Applying $(6)$, we get the contradiction $4 = 3$ and now the result follows since $g_{2}$ is injective. $\square$

\section{A Certain Family of Real Cyclotomic Fields} \label{Main}
For odd primes $p < q$ we denote by $\mathbb{Q}(\zeta_{pq})^{+}$ the real cyclotomic field of conductor $pq$. We write $h^{+}$ for its class number and $Cl^{+}$ for its ideal class group. We denote the corresponding $l$-parts by $h^{+}_{l}$ and $Cl^{+}_{l}$. As already stated in the Introduction, we ask that $p$ (or $q) \equiv 1~({\rm mod}~3)$ and that $f_{2} \neq f_{3}$, where  $f_{3} = p$ (or $q$) is the conductor of the real cyclic cubic subfield $K_{3}$ and $f_{2}$ is the conductor of any of the real quadratic subfields $K_{2}$ (hence $f_{2} = pq, p$ or $q$ ). Then $K_{6} = K_{2} K_{3}$ is a sextic subfield of  $\mathbb{Q}(\zeta_{pq})^{+}$ of conductor $pq$.

For the rest of this paper, we denote by $K_{6}$ the sextic field that we defined above, together with its subfields $K_{2}$ and $K_{3}$, with corresponding $3-$class numbers $h_{6,3}, h_{2,3}$ and $h_{3,3}$.

\begin{proposition}
\label{Prop1}
For the real cyclotomic field $\mathbb{Q}(\zeta_{pq})^{+}$ of $3-$class number exactly $9$, we assume that $p$ and $q$ satisfy the equivalence relations above. We assume further that $3 | h_{2,3}$.

If any one of the equations $(\ast)$ of Theorem ~\ref{Th2} has a solution, then the $3-$part of the ideal class group of $\mathbb{Q}(\zeta_{pq})^{+}$is cyclic, i.e. $Cl^{+}_{3} \cong \mathbb{Z}/9\mathbb{Z}$. 

\end{proposition}

\emph{Proof:} Consider the sextic field $K_{6}$ defined above. Regardless of whether $f_{6} = f_{2}$ or not, we see that from Theorem~\ref{Th3} we always have that $3 \mid w$ since we either have $w = 0$ (from Theorem~\ref{Th3}(i)) or $w = -2^{\nu + k + max \{ 1, \lambda \} - 2} h_{2}$ (from Theorem~\ref{Th3}(ii)) and $3 | h_{2}$ by assumption.

If $\chi_{3}(f_{3}) \neq 1$, Theorem~\ref{Th4}(ii) with $\lambda = \nu = 1$ holds, yielding $u^{2} - uv + v^{2} = 3^{k}$. For $k =1$, i.e. for $\chi_{6}$ nondecomposable, the equation $u^{2} - uv + v^{2} = 3$ gives $(u - v)^{2} = 3 - uv$ and $(u + v)^{2} = 3 (1 + uv)$. The two relations show that we must have $uv \leq 3$ and $1 + uv \geq 0$ therefore getting $-1 \leq uv \leq 3$. Assume $u=3$. Then $v = 1$ and substituting in our initial equation we get $3^2 - 3 \cdot 1 + 1^2 \neq 3$, contradiction. Similarly with $v$. Hence, the only possible values for $(u , v)$ are $\{ (1 , -1), (1 , 2), (-1 , 1), (-1 , -2), (2 , 1), (-2 , -1) \}$. Similar calculations for $k = 0$ give the pairs $\{ (0 , 1), (1 , 0), (1 , 1), (-1 , -1) \}$. We see that we never have $2 | u$ and $2 | v$ simultaneously. Theorem~\ref{Th1} now gives us that  $\langle - 1 \rangle N_{6/3}(U_{6}^{*}) = U_{3}$ and $N_{6/2}(U_{6}^{*}) \neq U_{2}$. Since we assumed that any one of the equations of Theorem ~\ref{Th2} has a solution, Theorem~\ref{Th2} now gives $[U_{6} : U_{6}^{*}] = 3$ and $N_{6/2}(U_{6}) = U_{2}$.

 If $\chi_{3}(f_{3}) = 1$, Theorem~\ref{Th4}(i) holds and $2$ divides both $u$ and $v$. From Theorem~\ref{Th1} we now have that  $\langle - 1 \rangle N_{6/3}(U_{6}^{*}) \neq U_{3}$ and $N_{6/2}(U_{6}^{*}) \neq U_{2}$. Regardless of whether any one of the equations $(\ast \ast)$ have a solution, Theorem~\ref{Th2(a)} gives $3 \parallel [U_{6} : U_{6}^{*}]$ and $N_{6/2}(U_{6}) = U_{2}$.
 
 We see that regardless of the value of $\chi_{3}(f_{3})$, both cases yield the same result: $3 \parallel [U_{6} : U_{6}^{*}]$ and $N_{6/2}(U_{6}) = U_{2}$.
 
 We know that $h_{6} = h_{2} h_{3} h_{R}$. Since $K_{3}/\mathbb{Q}$ is a ramified $3-$extension, we must have $h_{3} = 1$ (\cite{Wash1}, Theorem 10.4, p.186). Since $[U_{6} : U_{6}^{*}] = 3$ we have that $3 | h_{R}$ hence $9 | h_{6}$. Since $h^{+}_{3} = 9$, from \cite{Wash1} (Proposition 4.11, p. 39) we must have $h_{6,3} = 9$. 

Since there is no capitulation in $K_{6}/K_{2}$ (Theorem~\ref{3Capitulation} above), the index $[Cl_{6,3} : \mathfrak{t}(Cl_{2,3})] = 3^{1}$, where $\mathfrak{t} : Cl_{2,3} \mapsto Cl_{6,3}$ is the transfer of ideal classes mapping an ideal $\mathfrak{c}$ to $\mathfrak{c} O_{6}$. According to Theorem~\ref{Lem2} above, we must have that $\mathfrak{t}(Cl_{2,3}) = Cl_{6,3}^{3}$ therefore $3 = \mathfrak{t}(Cl_{2,3}) = Cl_{6,3}^{3}$ forcing $Cl_{6,3}$ to be cyclic. We thus have $Cl_{6,3} \cong \mathbb{Z}/9\mathbb{Z}$.

In the case $3 \nmid [\mathbb{Q}(\zeta_{pq}^{+}) : K_{6}]$ we immediately have that $Cl^{+}_{3} \cong \mathbb{Z}/9\mathbb{Z}$ from Note~\ref{Note} above. Assume otherwise. Consider any $3-$extension $K_{18}$ of $K_{6}$ contained  in ${Q}(\zeta_{pq})^{+}$. From Proposition 4.11 (p. 39 in \cite{Wash1}) applied twice (first to the extension $[K_{18} : K_{6}]$ and then to $[\mathbb{Q}(\zeta_{pq})^{+} : K_{18}]$), we must have that the $3-$class number of $K_{18}$ must also equal 9. But now, Proposition~\ref{Lem1} from above yields $r_{\mathfrak{t}} = r_{K_{6}} = 1$ hence we have capitulation. More specifically, from Proposition~\ref{Wash}, since we have a degree $3$ extension, the only classes that can capitulate are of order $3$ hence $ker(\mathfrak{t})_{K_{6} \mapsto K_{18}} \cong \mathbb{Z}/3\mathbb{Z}$. Since now we have $[Cl(K_{18})_{3} : \mathfrak{t}(Cl(K_{6})_{3})] = 3$, Proposition~\ref{Lem2} gives $Cl(K_{18})_{3}^{3} \cong \mathbb{Z}/3\mathbb{Z}$. Hence $Cl(K_{18})_{3}$ remains cyclic of order $9$. 

We apply the reasoning of the previous paragraph as many times as required to any extension $[\mathbb{Q}(\zeta_{pq})^{+} : K_{3^{i} \cdot 6}]$, $i \geq 0$. This procedure will eventually give us $Cl_{3}^{+} \cong \mathbb{Z}/9\mathbb{Z}$. $\square$   

The following is an immediate consequence of Proposition~\ref{Prop1}.

\begin{corollary}
\label{Corollary1}
The associated Hilbert $3-$class field tower of the fields of Proposition~\ref{Prop1} is finite of length 1. $\square$
\end{corollary}

\section{Examples} \label{Examples}

There are twelve sextic fields with incomplete data in the table of M$\ddot{a}$ki \cite{Maki1}, with two of them being of conductor $pq$ (\cite{Maki2}, pg. 607). However, very fast computations in PARI showed that none of these two fields has a real quadratic field with class number divisible by 3.  Therefore, we can say that from the data on Real Cyclic Sextic Fields in M$\ddot{a}$ki \cite{Maki1} and from Table 1 in Agathocleous \cite{Aga}, there are only four fields that fall into our category amongst all real cyclotomic fields of conductor $pq \leq 2021$. We present these fields below.

\subsection{The Real Cyclotomic Fields of conductor $3 \cdot 331$ and $7 \cdot 67$} Agathocleous \cite{Aga} calculated the $l$-part $h^{+}_{l}$ of $h^{+}$ for cyclotomic fields $\mathbb{Q}(\zeta_{pq})^{+}$ with $pq < 2000$ and for all odd primes $l < 10000$. There were eight cases of fields with $h^{+}_{l} > l$. Two of these cases, namely the real cyclotomic fields of conductor $7 \cdot 67$ and $3 \cdot 331$, fall in the category of fields that we study here. For both of these fields, their conductor $pq$ satisfies the equivalences of Proposition~\ref{Prop1} and their class number was found by \cite{Aga} to be exactly $9$. Furthermore, from very fast and simple computations in PARI/GP, their quadratic subfield $\mathbb{Q}(\sqrt{pq})$ was found to have $h_{2,3} = 3$. Finally, from the numerical data in M$\ddot{a}$ki \cite{Maki1} we have that for both of these fields , one of the equations of Theorem~\ref{Th2} has a solution $\xi_{B}$ in the sextic field $K_{6}$ defined in Proposition~\ref{Prop1}. Hence we conclude that $$Cl(\mathbb{Q}(\zeta_{3\cdot331})^{+})_{3} \cong Cl(\mathbb{Q}(\zeta_{7\cdot67})^{+})_{3} \cong \mathbb{Z}/9\mathbb{Z}.$$

\subsection{The Real Cyclotomic Fields of Conductor $7 \cdot 257$ and $3 \cdot 643$.} These fields satisfy all the conditions of Proposition ~\ref{Prop1} except from their $3-$class number, which is up to this point uknown. For the first one we have that $h_{2,3} = 3$ for $K_{2} = \sqrt{257} $ and for the second one $h_{2,3} = 3$ for $K_{2} = \sqrt{3 \cdot 643}$. Following the first steps outlined in the Proof of Proposition ~\ref{Prop1}, we see that the ideal class group of the sextic subfield $K_{6}$ for both of these cyclotomic fields is cyclic of order $9$. Also, since $3 \nmid [\mathbb{Q}(\zeta_{pq})^{+} : K_{6}]$, from Note ~\ref{Note}, we can state that $\mathbb{Z}/9\mathbb{Z} \cong \subseteq Cl(\mathbb{Q}(\zeta_{3 \cdot 643})^{+})_{3}$ and $\mathbb{Z}/9\mathbb{Z} \cong \subseteq Cl(\mathbb{Q}(\zeta_{7 \cdot 257})^{+})_{3},$ showing at the same time that the $3-$class number of these two real cyclotomic fields is divisible by 9.

We then follow the first two steps of the method outlined in \cite{Aga} and we find that the $3-$class number for both of these real cyclotomic fields cannot be greater than $9$, hence it must be equal to $9$. This result, together with Proposition~\ref{Prop1}, yield $$Cl(\mathbb{Q}(\zeta_{3\cdot643})^{+})_{3} \cong Cl(\mathbb{Q}(\zeta_{7\cdot257})^{+})_{3} \cong \mathbb{Z}/9\mathbb{Z}.$$

In the Appendix, we present the numerical data that we obtained from applying the first two steps of the method outlined in \cite{Aga}, for both of these fields. More specifically, we present the Frobenius polynomials obtained (by using PARI/GP) for the powers $M =3, 3^{2}$ and $3^{3}$ and the Gr\"obner Bases (computed in MATHEMATICA) for the ideals $J^{M}$. For a detailed explanation of what the Frobenius polynomials and the ideals $J^{M}$ stand for, the reader may refer to \cite{Aga}.

\begin{center}\textbf{Appendix}\end{center}

\textbf{The Field of Conductor $3 \cdot 643$}

Below we show the Frobenius polynomials obtained for $M =3, 3^{2}$ and $3^{3}$ for the pair of degrees $(d_{1},d_{2}) = (2,2)$:

\begin{center}
\tiny
\begin{tabular}{|c| p{11cm} |}
\hline
$l$ & The Frobenius Maps for $M$ = 3 \\
\hline \hline
$l_{1}=8437447$ & $f_{\Re_{1}} = (y^5 + y^3 + y^2 + 2y + 1)x + (2y^5 + y^4 + y^2 + 2y)$ \\
\hline
$l_{2} 11249929 = $ & $f_{\Re_{2}} = (y^5 + y^4 + 2y^3 + y^2 + y + 1)x + (2y^5 + 2y^4 + 2y^2 + y + 1)$ \\
\hline
$l_{3} 50624677 = $ & $f_{\Re_{3}} = 2y^3x + (y^5 + y^3 + y^2 + y)$ \\
\hline
$l_{4} = 78749497 $ & $f_{\Re_{4}} = (y^5 + y^4 + 2y^3 + y^2 + y)x + (y^4 + y^3 + 2y^2 + 2)$ \\
\hline
$l_{5} = 84374461$ & $f_{\Re_{5}} = (2y^5 + y + 1)x + 2y^5 $ \\
\hline
$l_{6} = 87186943$ &  $f_{\Re_{6}} = (y^5 + 1)x + (y^5 + y^4 + y^3 + 2y^2 + 2)$ \\
\hline
\hline
$l$ & The Frobenius Maps for $M = 3^{2}$\\
\hline \hline
$l_{1}=8437447$ & $f_{\Re_{1}} = (2y^5 + 5y^4 + y^3 + 3y^2 + 3y + 4)x + (3y^5 + 3y^4 + 4y^3 + 2y^2 + 5y + 1)$ \\
\hline
$l_{2} = 11249929 $ & $f_{\Re_{2}} = (6y^5 + 2y^4 + 3y^3 + 3y^2 + 3y + 1)x + (3y^5 + 3y^4 + y^3 + 6y^2 + 2y + 3)$ \\
\hline
$l_{3} = 50624677 $ & $f_{\Re_{3}} = (y^5 + y^4 + 5y^3 + 4y^2 + 7)x + (4y^5 + 7y^3 + y^2 + y + 5)$ \\
\hline
$l_{4} = 78749497 $ & $f_{\Re_{4}} = (7y^4 + 4y^3 + y^2 + 2y + 4)x + (y^5 + 2y^4 + 4y^3 + 7y + 4)$ \\
\hline
$l_{5} = 84374461$ & $f_{\Re_{5}} = (5y^5 + 3y^3 + 8y^2 + 4y + 7)x + (8y^5 + 4y^4 + 7y^3 + 5y^2 + 3)$ \\
\hline
$l_{6} = 87186943$ &  $f_{\Re_{6}} = (6y^5 + 6y^4 + 5y^3 + y^2 + 4y + 5)x + (y^5 + 4y^4 + 5y^3 + 6y^2 + 6y + 5)$ \\
\hline
\hline
$l$ & The Frobenius Maps for $M = 3^{3}$ \\
\hline \hline
$l_{1}=8437447$ & $f_{\Re_{1}} = (4y^5 + 21y^4 + 25y^3 + y^2 + 14y + 1)x + (11y^5 + 7y^4 + 21y^3 + 16y^2 + 20y + 21)$ \\
\hline
$l_{2} = 11249929$ & $f_{\Re_{2}} = (7y^5 + 13y^4 + 26y^3 + 25y^2 + y + 19)x + (5y^5 + 11y^4 + 9y^3 + 26y^2 + 7y + 13)$ \\
\hline
$l_{3} = 50624677$ & $f_{\Re_{3}} = (18y^5 + 18y^4 + 8y^3 + 12y^2 + 6)x + (19y^5 + 9y^4 + 10y^3 + 10y^2 + y + 24)$ \\
\hline
$l_{4} = 78749497 $ & $f_{\Re_{4}} = (y^5 + 16y^4 + 20y^3 + 13y^2 + y)x + (6y^5 + 10y^4 + 4y^3 + 17y^2 + 9y + 11) $ \\
\hline
$l_{5} = 84374461 $ & $f_{\Re_{5}} = (26y^5 + 12y^4 + 12y^2 + 16y + 13)x + (5y^5 + 6y^4 + 21y^3 + 24y^2 + 15y + 12)$ \\
\hline
$l_{6} = 87186943$ &  $f_{\Re_{6}} = (16y^5 + 9y^3 + 6y^2 + 6y + 22)x + (4y^5 + 25y^4 + y^3 + 17y^2 + 6y + 23)$ \\
\hline 
\end{tabular}
\end{center}

For $M = 3$, $J^{M} = (y^{2} - 1, y - x) \equiv ((y+1)(y-1), (y-1) - (x-1))$ with coefficients in $\mathbb{Z}/3\mathbb{Z}$. Analyzing the relations (as in \cite{Aga}), we see that $|I_{d}/J^{M}|$ = 3. 

For $M = 9$, $J^{M} = (-4+3y+y^2, 3+x-4y) \equiv ((y - 1)(y + 4), (x-1) - 4(y - 1))$ with coefficients in $\mathbb{Z}/9\mathbb{Z}$. Here we have that $|I_{d}/J^{M}| = 3^{2}$. Since $|I_{d}/J^{3}|$ is strictly smaller than  $|I_{d}/J^{3^{2}}|$ we need to continue as above with $M = 3^{3}$.

For $M= 27$, $J^{M} = (9-9y,4-3y-y^2, -3-x+4y)$ with coefficients in $\mathbb{Z}/27\mathbb{Z}$. We see here that $J^{3^{3}}$ is generated by the same polynomials as $J^{3^{2}}$ but it has the extra polynomial $9(y-1)$ which reduces the number of constants to 9 instead of 27. Therefore $|I_{d}/J^{3^{2}}|$ = $|I_{d}/J^{3^{3}}|$ = 9 and so, as expected, the orders of these quotients stabilize with $M = 3^{2}$. 

For the pair of degrees $(d_{1},d_{2})$ = (2,214) the Frobenius polynomials give exactly the same ideals $J^{M}$ as above and therefore we only need to consider the case for $(d_{1},d_{2}) = (2,2)$.

\textbf{The Field of Conductor $7 \cdot 257$}

This field has seven different pairs of degrees $(d_{1},d_{2})$. All give the same ideals $J^{M}$ and therefore we only consider the case for $(d_{1},d_{2}) = (2,2)$.  

\begin{center}
\tiny
\begin{tabular}{|c| p{11cm} |}
\hline
$l$ & The Frobenius Maps for $M$ = 3 \\
\hline \hline
$l_{1}=10491769$ & $f_{\Re_{1}} = x^5 + 2x^4 + x^2 + 2x$ \\
\hline
$l_{2} = 15737653$ & $f_{\Re_{2}} = (2y + 1)x^5 + yx^4 + 2x^3 + (2y + 1)x^2 + yx + 2$ \\
\hline
$l_{3} = 20983537$ & $f_{\Re_{3}} = (y + 2)x^5 + (y + 2)x^4 + (y + 2)x^3 + (y + 2)x^2 + (y + 2)x + (y + 2)$ \\
\hline
$l_{4} = 26229421$ & $f_{\Re_{4}} = (y + 2)x^5 + x^4 + 2yx^3 + (y + 2)x^2 + x + 2y$ \\
\hline
$l_{5} = 36721189$ & $f_{\Re_{5}} = 2yx^5 + 2x^4 + (y + 1)x^3 + 2yx^2 + 2x + (y + 1)$ \\
\hline
$l_{6} = 49835899$ &  $f_{\Re_{6}} = (2y + 2)x^5 + yx^4 + x^3 + (2y + 2)x^2 + yx + 1$ \\
\hline 
\hline
$l$ & The Frobenius Maps for $M = 3^{2}$\\
\hline \hline
$l_{1}=10491769$ & $f_{\Re_{1}} = x^5 + (6y + 5)x^4 + 6x^3 + x^2 + (6y + 5)x + 6$ \\
\hline
$l_{2} = 15737653$ & $f_{\Re_{2}} = (5y + 7)x^5 + (7y + 3)x^4 + (6y + 8)x^3 + (5y + 7)x^2 + (7y + 3)x + (6y + 8)$ \\
\hline
$l_{3} = 20983537$ & $f_{\Re_{3}} = (y + 5)x^5 + (4y + 2)x^4 + (4y + 2)x^3 + (y + 5)x^2 + (4y + 2)x + (4y + 2)$ \\
\hline
$l_{4} = 26229421$ & $f_{\Re_{4}} = (y + 5)x^5 + (6y + 1)x^4 + (2y + 3)x^3 + (y + 5)x^2 + (6y + 1)x + (2y + 3)$ \\
\hline
$l_{5} = 36721189$ & $f_{\Re_{5}} = 2yx^5 + (6y + 5)x^4 + (7y + 7)x^3 + 2yx^2 + (6y + 5)x + (7y + 7)$ \\
\hline
$l_{6} = 49835899$ &  $f_{\Re_{6}} = (8y + 2)x^5 + yx^4 + (3y + 4)x^3 + (8y + 2)x^2 + yx + (3y + 4)$ \\
\hline
\hline
$l$ & The Frobenius Maps for $M = 3^{3}$ \\
\hline \hline
$l_{1}=10491769$ & $f_{\Re_{1}} = (21y + 15)x^5 + (7y + 10)x^4 + 11yx^3 + (24y + 4)x^2 + (8y + 22)x + (16y + 24) $ \\
\hline
$l_{2} = 15737653$ & $f_{\Re_{2}} = (7y + 19)x^5 + (15y + 23)x^4 + (10y + 9)x^3 + (16y + 15)x^2 + (y + 25)x + (14y + 8)$ \\
\hline
$l_{3} = 20983537$ & $f_{\Re_{3}} = (9y + 22)x^5 + (23y + 2)x^4 + (25y + 22)x^3 + (19y + 19)x^2 + (8y + 9)x + (15y + 16)$ \\
\hline
$l_{4} = 26229421$ & $f_{\Re_{4}} = (y + 16)x^5 + (26y + 21)x^4 + (10y + 25)x^3 + 7x^2 + (7y + 25)x + (19y + 5)$ \\
\hline
$l_{5} = 36721189$ & $f_{\Re_{5}} = (25y + 23)x^5 + (6y + 1)x^4 + (5y + 13)x^3 + (4y + 4)x^2 + 22x + (2y + 3)$ \\
\hline
$l_{6} = 49835899$ &  $f_{\Re_{6}} = (12y + 2)x^5 + (15y + 25)x^4 + (15y + 12)x^3 + (14y + 18)x^2 + (4y + 20)x + (15y + 10)$ \\
\hline
\end{tabular}
\end{center}

For $M = 3$, $J^{M} = (-1+y^2, 1-x) \equiv ((y+1)(y-1), (x-1))$ with coefficients in $\mathbb{Z}/3\mathbb{Z}$, which gives $|I_{d}/J^{M}|$ = 3. 

For $M = 9$, $J^{M} = (-1+y^2, 2+x-3y) \equiv ((y-1)(y+1), (x-1) - 3(y-1)$ with coefficients in $\mathbb{Z}/9\mathbb{Z}$, which gives $|I_{d}/J^{M}| = 3^{2}$.

For $M= 27$, $J^{M} = (-9+9y, -1+y^2, 2+x-3y)$ with coefficients in $\mathbb{Z}/27\mathbb{Z}$. We see here that $J^{3^{3}}$ is generated by the same polynomials as $J^{3^{2}}$ but it has the extra polynomial $9(y-1)$ which reduces the number of constants to 9 instead of 27. Therefore $|I_{d}/J^{3^{2}}|$ = $|I_{d}/J^{3^{3}}|$ = 9 and the orders of these quotients stabilize with $M = 3^{2}$.

\vspace{10mm}

\begin{flushleft}Eleni Agathocleous \\ University of Nicosia, Cyprus \\ email: agathocleous.e@unic.ac.cy \end{flushleft}


\begin{thebibliography}{99}
\bibitem{Aga} E. Agathocleous: On the Class Numbers of Real Cyclotomic Fields of Conductor pq, {\em Acta Arithmetica} \textbf{165(3)} (2014), 257-277.
\bibitem{B1} C. Batut, K. Belabas, D. Bernardi, H. Cohen, M. Olivier: User's Guide to PARI-GP (version 2.3.0),
Universit\'{e} Bordeaux I, Bordeaux 2000.
\bibitem{Gras} M-N. Gras: M\'ethodes et algorithmes pour le calcul num\'erique du nombre de classes et de unit\'es des extensions cubiques cycliques de Q, {\em Journal f\"ur die reine und angewandte Mathematik} \textbf{277} (1975), 89-116.
\bibitem{Maki2}  V. Ennola, S. M$\ddot{a}$ki, R. Turunen:  On Real Cyclic Sextic Fields, {\em Mathematics of Computation} \textbf{45(172)} (1985), 591-611.
\bibitem{Lemm2} F. Lemmermeyer: Ideal class groups of cyclotomic number fields II, {\em Acta Arithmetica} \textbf{84(1)} (1998), 59-70.
\bibitem{Maki1} S. M$\ddot{a}$ki, {\em The Determination of Units in Real Cyclic Sextic Fields}, Springer-Verlag, Berlin Heidenberg GmbH, 1980.
\bibitem{Wash2} R. Schoof, L.C. Washington: Visibility of Ideal Classes, {\em Journal of Number Theory} \textbf{130(12)} (2010), 2715-2731.
\bibitem{Wash1} L.C. Washington, {\em Introduction to Cyclotomic Fields; second edition}, Springer-Verlag, Berlin Heidenberg New York, 1997.
\bibitem{W2} Wolfram Research, Inc., Mathematica, Version 10.0, Champaign, IL (2014).
\end{thebibliography}
\end{document}